\documentclass[12pt]{amsart}
\usepackage{graphicx}
\usepackage{amssymb}
\usepackage{epstopdf}
 \usepackage{latexsym}
\usepackage{amsmath}
\usepackage{amsfonts} 
\usepackage{amsthm}

\newtheorem{Proposition}{Proposition}
\newtheorem{Theorem}[Proposition]{Theorem}
\newtheorem{Lemma}[Proposition]{Lemma}

\newtheorem{Corollary}[Proposition]{Corollary}

\def\XXint#1#2#3{{\setbox0=\hbox{$#1{#2#3}{\int}$}
\vcenter{\hbox{$#2#3$}}\kern-.5\wd0}}

\def\e{\epsilon}

    \def\sqr#1#2{{\vcenter{\vbox{\hrule height .#2pt
                             \hbox{\vrule width .#2pt height#1pt \kern#1pt
                                   \vrule width .#2pt}
                             \hrule height .#2pt}}}}

     \def\CC{\mathbb{C}}

    \def\e{\epsilon}

\def\be{\begin{equation}}
\def\ee{\end{equation}}

\begin{document} 

\title[First return map ]{The first return map for planar vector fields with nilpotent linear part with a center or a focus}
\author{Rodica D. Costin}

\maketitle


\begin{abstract}

The return map for planar vector fields with nilpotent linear part (having a center or a focus and under an assumption generically satisfied) is found as a convergent power series whose terms can be calculated iteratively. The first nontrivial coefficient is the value of an Abelian integral, and the following ones are explicitly given as iterated integrals built with algebraic functions.

\end{abstract}


\section{Introduction}

The study of planar vector fields has been the subject of intense investigation, due to their importance in applications, and in connection to Hilbert's 16th Problem \cite{Rouss_3}. Significant progress has been made in the geometric theory of these fields, as well as in bifurcation theory, normal forms, foliations,  and the study of Abelian integrals (see the recent book \cite{Dum_Book}).

First return maps have been studied in relation to existence of closed orbits; more generally, return maps are important in a large array of applications (see \cite{Guckenheimer} and references therein) and also in logic, in connection to o-minimality \cite{Miller}.

A fundamental result regarding the asymptotic form of return maps states that if the singular points of a $C^\infty$ vector field are algebraically isolated, there exists a semitransversal arc such that the return map admits an asymptotic expansion is positive powers of $x$ and logs (with the first term linear), or has its principal part a finite composition of powers and exponentials \cite{Ilyashenko},\,\cite{Moussu}.

In the case when the linear part of the vector field has non-zero eigenvalues there is a good understanding of the return map \cite{Brudnyi},\,\cite{Dumortier},\,\cite{Rouss_1},\,\cite{Rouss_2},\,\cite{Francoise_Lienard_2}, \cite{Francoise_Lienard_1}, as well as for perturbations of Hamiltonians \cite{Francoise_Iterated_int},\,\cite{Francoise_Suc_der}, or for perturbations of integrable systems \cite{Gavrilov}.
In the general setting, however, there are few results available \cite{Berezovskaya},\,\cite{RDC_Ret1},\,\cite{Medvedeva_Lt},\,\cite{Medvedeva_St}. 

The present paper establishes an iterative procedure for calculating the return map as an integer power series for generic vector fields with nilpotent linear part in the case of a center or a focus. The first few coefficients in these series are explicitly given. Using an algorithm given here these coefficients can be found explicitly up to any (finite) order in terms of iterated integrals involving algebraic functions.

A fundamental result concerning planar vector fields analytic near a stationary point $(0,0)$ states that if the linear part is not zero, having both eigenvalues zero, then such fields have the normal form near $(0,0)$
\be\label{fnf}
\dot{x}=y,\ \dot{y}=a(x)+yb(x)
\ee
with $a(x), b(x)$ analytic at $0$, with $a$ (respectively, $b$) having a zero of order at least two (respectively, one) at $x=0$ \cite{I_Y_Book}.

Furthermore, if the origin is a center or a focus, and if the linear part of this system is not radial (i.e. does not have the form $\dot{x}=\lambda x,\, \dot{y}=\lambda y$), then after an analytic change of variables (\ref{fnf}) can be written as \cite{Loray}
\be\label{cent_foc}
\dot{z}\,=\, -\,w\,f(z)\,+\,z^{l+1}\,g(z),\ \ \ \ \ \ \dot{w}\,=\,k\,z^{2k-1}\,f(z)+k\,w\,z^l\,g(z)
\ee
where 
\be\label{as1}
1\leq k\leq l+1\ \ \ \  \ \ \ {\mbox{and}} \ \ \ \  \ \ \ f(0)\ne 0
\ee
The present paper studies the first return map for (\ref{cent_foc}) under the supplementary {\bf{assumption}}:
\be\label{as2}
 k\ne l+1
 \ee

 \section{Main result}

The main result is the following:

\begin{Theorem}\label{MainProp} 

\ 

Consider the system (\ref{cent_foc}) with $f(z), g(z)$ analytic at $0$ and satisfying (\ref{as1}),\,(\ref{as2}), and its transversal $\mathcal{T}:\ f(z)\, w=z^{l+1}\, g(z)$ ($z>0$, small).

Then the first return to $\mathcal{T}$ of the solution with $z(0)=\e$ (with $\e$ small enough) is at the point with $z=Z(\e)$ analytic in $\e$, with 
$$Z(\e)=\e+\sum_{n\geq l-k+2}\, \e^n\, Z_n$$
where $Z_n$ can be calculated iteratively in terms of iterated integrals involving algebraic functions.

\end{Theorem}

The proof of {Theorem}\,\ref{MainProp} shows that $\mathcal{T}$ is indeed transversal.
Also, the proof contains, and it is built around, an algorithm for calculating iteratively the coefficients $Z_n$; in particular, the first two of them are calculated: (\ref{zretp}), (\ref{zret1}) with the notations (\ref{not_p}),\,(\ref{def_F}),\,(\ref{not_PHI}),\,(\ref{valB0}),\,(\ref{tetp}).

Necessary conditions for a system to have closed orbits (thus a center) are discussed in \S\ref{center}.

 \section{Proof of Theorem\,\ref{MainProp}}

\subsection{Normalization}\label{normalization}

Eliminating the time in (\ref{cent_foc}) we obtain
\be\label{eqwz}
\frac{dw}{dz}=\frac{kz^{2k-1}f(z)+kwz^lg(z)}{-wf(z)+z^{l+1}g(z)}
\ee

Solutions of (\ref{cent_foc}) provide smooth parametrizations for solutions of (\ref{eqwz}). There are points where the graph of a solution $w(z)$ of (\ref{eqwz}) has vertical tangents; at such a point the parametrization (\ref{cent_foc}) prescribes that $w(z)$ is continued with another solution through the same point, and having a vertical tangent as well.

It is convenient to straighten the curve $\mathcal{T}$ using a substitution: let
\be\label{subwu}
w=z^{l+1}\, F(z)+w_1
\ee
with
\be\label{def_F}
F(z)= \frac{g(z)}{f(z)}
\ee
(note that $F(z)$ is analytic at $0$ since $f(0)\ne 0$). Equation (\ref{eqwz}) becomes
\be\label{eqw1}
 \frac{d(w_1^2)}{dz}= -2kz^{2k-1}A(z)-w_1z^{p+k-1} \, B(z)
\ee
with the notations 
\be\label{not_p}
p=l-k+1,\ \ \ p\geq 1
\ee
\be\label{notA}
A(z)\,=\,1+z^{2p}\, F^2(z)\, \equiv\, A_F(z)
 \ee
 \be\label{notB}
B(z)\,=\,2\,z\frac{dF}{dz}+2(p+2k)F(z) \, \equiv\, B_F(z)
 \ee

The problem translates into the study of the first return to the positive $z$-axis of solutions of (\ref{eqw1}) which start close enough to the origin, i.e. which satisfy the initial condition $w(\epsilon)=0$ with  $\epsilon>0$ small enough.

It is convenient to introduce the small parameter $\epsilon$ in the equation. With the notation
\be\label{notxy}
z=\epsilon x,\ w_1=\epsilon^ky
\ee
equation (\ref{eqw1}) becomes
\be\label{eqy}
\frac{d(y^2)}{dx}=-2k\,x^{2k-1}A(\e x)-\e^{p}y\,x^{p+k-1}\, B(\e x)
\ee

While the initial condition $w_1(\e)=0$ becomes $y(1)=0$, it is useful to study solutions with the more general initial condition $y(\eta)=0$ with $\eta$ in a neighborhood of $1$.


It is assumed that 
\be\label{eta}
\eta\in(1-\eta_0,1+\eta_0)\ \ \ \ {\mbox{for\ some\ fixed\ }} \eta_0\in(0,1)
\ee

\subsection{General behavior of solutions of (\ref{eqy}).} 

In parametrized form equation (\ref{eqy}) is 
\be\label{eqxyt}
\dot{x}=-2y,\ \dot{y}=2k\,x^{2k-1}\,A(\e x)+\e^{p}y\,x^{p+k-1}\, B(\e x)
\ee

Consider the solution of (\ref{eqxyt}) with the initial condition $x(0)=\eta,\, y(0)=0$ for some $\eta>0$. Let $(x_0,y_0)$ be the solution of the system (\ref{eqxyt}) for $\e=0$:
$$\dot{x_0}=-2y_0,\ \dot{y_0}=2k\,x^{2k-1},\ \ \ x_0(0)=\eta,\, y_0(0)=0$$
This solution satisfies
\be\label{sol0}
y_0^2+x_0^{2k}=\eta^{2k}
\ee
and clearly 
 \be\label{sol_close}
 x(t)=x_0(t)+O(\e),\ \ \ \ \ y(t)=y_0(t)+O(\e)\ \ \ \ \ \  \ \ \ \ (\e\to 0)
 \ee
  for any finite interval of time.

\subsection{Positive solutions of (\ref{eqy}) for $x>0$.}\label{QI}

Lemma\,\ref{L1} shows that there exists a unique solution of (\ref{eqy}) so that $y(\eta)=0$ and $y\geq 0$ for $\e=0$. It also shows that this solution is defined for $x\in[0,\eta]$ and establishes an iterative procedure for calculating its power series expansion in $\e$.

With the substitution $y=u^{1/2}$ equation (\ref{eqy}) becomes
\be\label{equ}
\frac{du}{dx}=-2k\,x^{2k-1}A(\e x)-\e^{p}u^{1/2}\,x^{p+k-1}\, B(\e x)
\ee

\begin{Lemma}\label{L1}

There exists $\e_0>0$ so that so that the following holds.

Let $\eta$ with (\ref{eta}). For any $\e$ with $|\e |<\e_0$ equation (\ref{equ}) with the condition $u(\eta)=0$ has a unique solution $u=u(x;\e,\eta)$ for $x\in [0,\eta]$. We have $u(x;0,\eta)>0$ for $x\in[0,\eta)$,  and $u(x;\e,\eta)$ is analytic in $\e$ and $\eta$.

\end{Lemma}

{\em{Proof of Lemma\,\ref{L1}.}}

Denote
\begin{multline}\label{defpsi}
{\tilde{P}}(x;\eta,\e)=\frac{1}{\eta-x}\,\int_\eta^x\, \left[ -2kt^{2k-1}\, A(\e t)\right]\, dt\\
=\frac{\eta^{2k-1}}{1-x/\eta}\,\int_{x/\eta}^1\, 2k \sigma^{2k-1}\, A\left(\e \eta \sigma\right)\, d\sigma
\end{multline}

{\em{Heuristics.}} A particular vector field was studied in \cite{RDC_Ret1}, unveiling the main ideas involved. In the present general case a few heuristic considerations are the following.

Looking for solutions of (\ref{eqy}) with $y(\eta)=0$ we obtain that
$$y(x)\sim c(\eta-x)^{1/2} (1+o(1))\ \ (x\to\eta)\ \ \ {\mbox{with\ }} c^2=2k\eta^{2k-1}A(\e\eta)$$
Note that $A(\e\eta)=1+O(\e)>0$ for $\e$ small enough, in view of (\ref{notA}).
Therefore 
\be\label{heury}
u=y^2=(\eta-x)\tilde{P}(x;\eta,\e) +\Delta(x)
\ee
{with} $\Delta=o(\eta-x)$ as $x\to\eta$.
 In fact, substitution of (\ref{heury}) in (\ref{equ}) gives 
 \be\label{heury2}
 \Delta\sim {\rm{const}}\ (\eta-x)^{3/2}
 \ee

These observations motivate the following substitutions.

Denote 
$$\xi=\sqrt{1-\frac{x}{\eta}}$$
(with the usual branch of the square root for $x/\eta<1$).
Let
\be\label{def_delta}
\delta=\e\eta
\ee
Note that (\ref{defpsi}) can be written
$\tilde{P}(x;\eta,\e)\equiv \eta^{2k-1}P(\xi;\delta)$
where
\begin{multline}\label{form_psi}
P(\xi;\delta)=\xi^{-2}\,\int_{1-\xi^2}^1\, 2k\,\sigma^{2k-1}\, A\left(\delta\sigma\right)\, d\sigma\\
=2k\int_0^1(1-\xi^2s)^{2k-1}A\left(\delta(1-\xi^2s)\right)\,ds
\end{multline}

With the substitution
\be\label{subuv}
u \left( x \right) =(\eta-x)\eta^{2k-1}\left[P(\xi;\delta) + v\left( \xi;\delta \right)\right] 
\ee
equation (\ref{equ}) becomes
\be\label{eqv}
\xi\,{\frac {dv}{d\xi}} +2\,v  =
2\delta^p \,{\xi}(1-\xi^2)^{p+k-1}\,(P+v)^{1/2} \, B\left(\delta(1-\xi^2)\right)
\ee

Note that $y(\eta)=0$ implies $u(\eta)=0$ and then necessarily $v(0)=0$ by (\ref{defpsi}), (\ref{heury}), (\ref{heury2}).

Lemma\,\ref{L1} follows with $\e_0=\delta_0/(1-\eta_0)$ if we show the following:

\begin{Lemma}\label{Lserv}

These exists $\delta_0>0$ so that for any $\delta$ with $|\delta|<\delta_0$ equation (\ref{eqv}) has a unique solution $v=v(\xi;\delta)$ so that $v(0)=0$ and this solution is defined for $\xi\in[0,1]$.

Moreover, $v(\xi;\cdot)$ is analytic for $\delta\in\CC$ with $|\delta |<\delta_0$. 

The terms of its power series 
\be\label{servxi}
 v(\xi;\delta)=\sum_{n\geq p}\delta^nv_n(\xi)
 \ee
 can be calculated recursively. In particular, the first terms are as follows.
 
 Using the notations\footnote{$\Phi_{p,k}(\xi)$ is expressible in terms of the incomplete beta function.}
  \be\label{not_PHI}
  \Phi_{p,k}(\xi)= \xi^{-2}\int_{1-\xi^2}^1s^{p+k-1}(1-s^{2k})^{1/2}\, ds
  \ee

  \be\label{not_PSI}
  \Psi_k(\xi)\,=\,\xi^{-2}\int_{1-\xi^2}^1ds\, \frac{s^{k}}{(1-s^{2k})^{1/2}}\,\int_s^1d\sigma\, \sigma^k(1-\sigma^{2k})^{1/2}
  \ee
and 
\be\label{valB0}
B_0\,=\,2(p+2k)\, F(0)\,\equiv\, B_{F;0},\ \ \ B_1\,=\,2\,(p+2k+1)F'(0)\,\equiv\, B_{F;1}
\ee
we have:
 \be\label{vp}
  v_p(\xi)=\, B_0 \Phi_{p,k}(\xi)
  \ee
 and, moreover, for $p\geq 2$ we have
 \be\label{vpnext}
 v(\xi;\delta)=\delta^{p}\, v_p(\xi)
 +\delta^{p+1}\, B_1\, \Phi_{p+k+1,k}(\xi)
 +O\left(\delta^{p+2}\right) 
 \ee
 while for $p=1$ we have
 \be\label{vp1}
 v(\xi;\delta)=\delta\, v_1(\xi)
 +\delta^{2}\,\left[ \,B_1\, \Phi_{2,k}(\xi)+\,\frac{1}{2}\,{B_0^2}\,\Psi_k(\xi)\,\right]
 +O\left(\delta^{3}\right) 
 \ee

 \end{Lemma}

{\em{Proof of Lemma\,\ref{Lserv}.}} 

Multiplying (\ref{eqv}) by $\xi$ and integrating we obtain:
$$v(\xi;\delta)={\xi}^{-2}\left[ C+2\delta^p \,\int_0^\xi t^2(1-t^2)^{p+k-1}\left[P(t;\delta)-v(t;\delta)\right]^{1/2}\, B\left(\delta(1-t^2)\right)\, dt\right]$$
Since $v(0)=0$ then the constant $C$ must vanish. It follows that $v$ is a fixed point ($v=\mathcal{J}[v]$) for the operator
\begin{multline}\label{defJ}
 \mathcal{J}[v]\,\left( \xi;\delta \right) =2\delta^p \,{\xi}^{-2}\int_0^\xi t^2(1-t^2)^{p+k-1}\left[P(t;\delta)+v(t;\delta)\right]^{1/2}\, B\left(\delta(1-t^2)\right)\, dt\\
 =2\delta^p \,{\xi}\int_0^1 s^2(1-\xi^2s^2)^{p+k-1}\left[P(\xi s;\delta)+v(\xi s;\delta)\right]^{1/2}\, B\left(\delta(1-\xi^2s^2)\right)\, ds\
 \end{multline}

{\em{Estimates for $P(t;\delta)$:}}

Let $r>0$ be small enough so that $F$ is analytic in a neighborhood of $|z|\leq r$, and so that $|A(z)-1|\leq c_0$ for all $|z|\leq r$ for some $c_0$ with $0<c_0<1/(2k)$ (see (\ref{notA})).

  Then from (\ref{form_psi}), for $\xi\in[0,1]$ and $|\delta|\leq r$ we have the upper estimate
 \be\label{upestP}
  |P(\xi;\delta)|\,\leq 2k(1+c_0)
  \ee
 and the lower estimate
 \begin{multline}
 |P(\xi;\delta)|\geq 2k\int_0^1(1-\xi^2s)^{2k-1}\,ds-\big|2k\int_0^1(1-\xi^2s)^{2k-1}\left[A\left(\delta(1-\xi^2s)\right)-1\right]\,ds\big|\\
 \geq 2k\int_0^1(1-s)^{2k-1}\,ds-2k\int_0^1(1-\xi^2s)^{2k-1}\big|A\left(\delta(1-\xi^2s)\right)-1\big|\,ds\\
 \geq 1-2kc_0>0
 \end{multline}

{\em{The operator $ \mathcal{J}$ is contractive:}}

Let $M=\sup_{|z|\leq r}|B(z)|$.

Let $\mu$ be a number with $0<\mu<1-2kc_0$.

 Let $\delta_0>0$ be small enough, so that 
 \be\label{choosed0}
  \delta_0^p\,\frac{2}{3}\, [2k(1+c_0)+\mu]^{1/2}M<\mu,\ \ \ \delta_0^p\frac{M}{3(1-2kc_0-\mu)^{1/2}}<1
  \ee
 
Let $\mathcal{B}$ be the Banach space of functions $f(\xi;\delta)$ continuous for $\xi\in [0,1]$ and analytic on the (complex) disk $|\delta|<\delta_0$, continuous on $|\delta|\leq\delta_0$, with the norm 
$$\|f\|=\sup_{\xi\in[0,1]}\,\sup_{|\delta|\leq\delta_0}\, |f(\xi;\delta)|$$

Let $\mathcal{B}_\mu$ be the ball $\mathcal{B}_\mu=\{ f\in \mathcal{B}\, ;\, \|f\|\leq \mu\}$. 

The operator $\mathcal{J}$ defined by (\ref{defJ}) is defined on $\mathcal{B}_\mu$. Indeed, for $f\in \mathcal{B}_\mu$ we have
\be\label{plus}
|P(t;\delta)+f(t;\delta)|\geq |P(t;\delta)|-| f(t;\delta)|\geq 1-2kc_0-\mu>0
\ee
therefore, since $f$ and $P$ are analytic in $\delta$, then so is $\left( P+f   \right) ^{1/2}$, and therefore so is $\mathcal{J}[f]$. 

Also, $\mathcal{J}\mathcal{B}_\mu\subset\mathcal{B}_\mu$ since for $f\in \mathcal{B}_\mu$, using (\ref{defJ}), (\ref{upestP}), (\ref{choosed0}) we have
$$\big| \mathcal{J}f\,(\xi;\delta) \big|\leq  |\delta|^p\,\frac{2}{3}\, [2k(1+c_0)+\mu]^{1/2}M<\mu$$

Moreover, the operator $\mathcal{J}$ is a contraction on $\mathcal{B}_\mu$. Indeed, using the estimate
$$\big| (P+f_1)^{1/2}-(P+f_2)^{1/2}\big|\leq |f_1-f_2|\, \frac{1}{2}\,\sup_{f\in \mathcal{B}_\mu}|P+f|^{-1/2}\leq \frac{|f_1-f_2|}{2(1-2kc_0-\mu)^{1/2}}$$
(by (\ref{plus})) we obtain
$$\big| \mathcal{J}f_1-\mathcal{J}f_2 \big|\leq c\|f_1-f_2\|\ \ \ {\mbox{with}}\ c=\delta_0^p\frac{M}{3(1-2kc_0-\mu)^{1/2}}<1$$
(by (\ref{upestP})).

Therefore the operator $\mathcal{J}$ has a unique fixed point in $\mathcal{B}_\mu$, which is the solution $v(\xi;\delta)$, analytic for $\delta\in\CC$ with $|\delta |<\delta_0$.

{\em{A recursive algorithm for calculating the power series in $\e$:}}

To obtain the power series (\ref{servxi}) substitute an expansion $ v(\xi;\delta)=+\sum_{n\geq 0}\delta^nv_n(\xi)$ in (\ref{eqv}). It follows that for $n<p$ we have $\xi v_n'+2 v_n=0$ with $v_n(0)=0$, therefore $v_n(\xi)\equiv 0$ for $n<p$.

Using (\ref{form_psi}), (\ref{notA}) we obtain for $P(\xi;\delta)$ a power series in $\delta$, with coefficients polynomials in $\xi$:
\be\label{serp}
P(\xi;\delta)=P_0(\xi)+\sum_{m\geq 0}\delta^{2p+m}P_{2p+m}(\xi)
\ee
where
\be\label{formP0}
{P_0} \left( \xi \right) ={\frac { 1-  \left( 1-{\xi}^{2} \right) ^{2k}}{{\xi}^{2}}}
\ee
and 
\be\label{formPm}
 {P_{2p+m}} \left( \xi \right) =\, F_{2,m}\, 
\frac{2k}{2p+2k+m}\, \frac {  1-   \left( 1-{\xi}^{2} \right) ^{2\,k+2\,p+m}
  }{\xi^{2} }
  \ee
with the notation 
\be\label{serFsq}
F^2(z)\equiv\sum_{m\geq 0} F_{2,m}z^m
\ee
 and in particular
\be\label{formF20}
F_{2,0}=F(0)^2
\ee

Substitution of (\ref{servxi}), (\ref{serp}), followed by expansion in power series in $\delta$ give
\be\label{stareq}
(P(\xi;\delta)+v(\xi;\delta))^{1/2} \, B\left(\delta(1-\xi^2)\right)\equiv\sum_{n\geq 0}\delta^nR_n(\xi)
\ee
where $R_n=R_n[v_p,\ldots,v_n,\xi]$.

From (\ref{eqv}) and (\ref{stareq}) we obtain the recursive system
$$\xi\,{\frac {dv_n }{d\xi}}+2\,v_n  =2{\xi}(1-\xi^2)^{p+k-1}\, R_{n-p},\ \ \ \ \  \ n\geq p$$
with the only solution with $v_n(0)=0$ given recursively by
\be\label{recvn}
v_n(\xi)=2\xi^{-2}\int_0^\xi\, {t^2}(1-t^2)^{p+k-1}\,R_{n-p}(t)\, dt
\ee

{\em{The first terms:}}

To calculate the first few $R_n$ note that
\begin{multline}\label{firtwo}
(P(\xi;\delta)+v(\xi;\delta))^{1/2} \, B\left(\delta(1-\xi^2)\right)\\=B_0P_0(\xi)^{1/2}+\delta B_1P_0(\xi)^{1/2} (1-\xi^2)+\delta^p\frac{B_0v_p(\xi)}{2P_0(\xi)^{1/2}}+O\left(\delta^2\right)
\end{multline}
where $B_0,B_1$ are coefficients in the expansion $B(z)=B_0+zB_1+O(z^2)$, and in view of (\ref{notB}), they are (\ref{valB0}).

In particular
 \begin{multline}
 v_p(\xi)=2B_0\,\xi^{-2}\int_0^\xi t^2(1-t^2)^{p+k-1}P_0(t)^{1/2}\, dt\\
 =B_0\xi^{-2}\int_{1-\xi^2}^1s^{p+k-1}(1-s^{2k})^{1/2}\, ds
 \end{multline}
which equals (\ref{vp}).

 {{For $p\geq 2$}} relation (\ref{firtwo}) is
\be\label{ptwo}
(P+v)^{1/2} \, B\left(\delta(1-\xi^2)\right)\\=B_0P_0(\xi)^{1/2}+\delta B_1P_0(\xi)^{1/2} (1-\xi^2)+O\left(\delta^2\right)
\ee
and using (\ref{firtwo}) we obtain (\ref{vpnext}).

{{For $p=1$}} we have $O(\delta^{2p})=O(\delta^{p+1})$ and there is one more term in the second nontrivial coefficient of $v(\xi;\delta)$. The calculation is straightforward: 
using relation (\ref{firtwo}) for $p=1$ and (\ref{recvn}) we obtain (\ref{vp1}).

\qed

\

The following {Corollary} gathers the conclusions of the present section. Many quantities depend on the function $F$, see (\ref{def_F}), (\ref{notA}), (\ref{notB}), (\ref{form_psi}), and we add the subscript $F$ for them:

\begin{Corollary}\label{Coro1}
There exists $\e_0>0$ so that for any $\eta$ with (\ref{eta}) and $\e$ with $|\e |<\e_0$ equation (\ref{eqy}) has a unique solution $y(x)$ on $[0,\eta]$ satisfying $y(\eta)=0$ and $y>0$ on $[0,\eta)$ for $\e=0$. 

Moreover, this solution has the form
\be\label{formphi}
y(x)=\phi_F(x;\e,\eta)\equiv(\eta-x)^{1/2}\eta^{k-1/2}\left[P_F\left( \sqrt{1-\frac{x}{\eta}} ;\e{\eta} \right) + v_F\left(  \sqrt{1-\frac{x}{\eta}} ;\e{\eta} \right)\right] ^{1/2} \ee
with $P_F(\xi;\delta)\equiv P(\xi;\delta)$ given by (\ref{form_psi}) and $v_F(\xi;\delta)$ solution of (\ref{eqv}) with $v_F(0;\delta)=0$ and $v_F(\xi;0)=0$. 

The map $(\e,\eta)\mapsto v_F(\xi;\e{\eta})$ is analytic for $|\e|<\e_0$ and $\eta$ as in (\ref{eta}).
\end{Corollary}
\qed

\subsection{Solutions of (\ref{eqy}) in other quadrants and matching}\label{QA}

\subsubsection{Solutions in the four quadrants}

{Corollary}\,\ref{Coro1} gives an expression for the solution $y(x)$ of (\ref{eqy}) for $x>0$ and $y>0$. Solutions in the other quadrants are found in the following way.

Let $\e$ with $|\e|<\e_0$ with $\e_0$ given by Lemma\,\ref{L1}.

{\bf{(i)}} Let $y_1(x)=\phi_F(x;\e,\eta)$ the solution (\ref{formphi}) of (\ref{eqy}), defined for $x\in[0,\eta]$, with $y_1(\eta)=0$. We have $y_1(x)=(\eta^{2k}-x^{2k})^{1/2}+O(\e)$ in view of (\ref{sol0}), (\ref{sol_close}), and we will refer to $y_1$ as a "solution in the first quadrant".

Solutions "in the other quadrants" are obtained as follows.

For a function $F(z)$ denote by $J_iF$ the following functions:
\begin{multline}\label{def_Ji}
(J_2F)(x)=(-1)^{p+k}F(-x),\ (J_3F)(x,y)=(-1)^{p+k-1}F(-x),\\ (J_4F)(x)=-F(x)
\end{multline}

Note that $J_2J_3=J_4$ and $A_{J_2F}=A_F$.

{\bf{(ii)}} It is easy to check that the function $y_2(x)=\phi_{J_2{F}}(-x;\e,\eta)$ is a solution of (\ref{eqy}). It is obviously defined for $x\in[-\eta,0]$; we have $y_2(-\eta)=0$ and $y_2(x)=(\eta^{2k}-x^{2k})^{1/2}+O(\e)$.

{\bf{(iii)}}  Similarly, the function $y_3(x)=-\phi_{J_3{F}}(-x;\e,\eta)$  is a solution of (\ref{eqy}) defined for $x\in[-\eta,0]$. We have $y_3(-\eta)=0$ and $y_3(x)=-(\eta^{2k}-x^{2k})^{1/2}+O(\e)$.

{\bf{(iv)}} The function $y_4(x)=-{\phi_{J_4{F}}(x;\e,\eta)}(x;\e,\eta)$ is a solution of (\ref{eqy}) defined for $x\in[0,\eta]$; we have $y_4(\eta)=0$ and $y_4(x)=-(\eta^{2k}-x^{2k})^{1/2}+O(\e)$.

\subsubsection{Matching at the positive $y$-axis}

Let $\eta,\tilde{\eta}$ satisfying (\ref{eta}) and let $y_1(x)=\phi_F(x;\e,\eta)$ solution as in {\bf{(i)}}, for $x\in[0,\eta]$ and $\tilde{y_2}(x)=\phi_{J_2F}(-x;\e,\tilde{\eta})$ solution as in {\bf{(ii)}}, for $x\in[-\tilde{\eta},0]$. 

The following Lemma finds $\tilde{\eta}$ so that $y_1(0)=\tilde{y_2}(0)$, therefore so that $y_1$ is the continuation of $\tilde{y_2}$:

\begin{Lemma}\label{Lyp}

Let $|\e|<\e_0$. Let $\eta$ so that $|\eta-1|\leq c_1|\e|$ with $c_1$ small enough so that $c_1\e_0< \eta_0$.

There exists a unique $\tilde{\eta}=\eta+O(\e)$ so that 
\be\label{eqphi1}
\phi_F(0;\e,\eta)=\phi_{J_2{F}}(0;\e,\tilde{\eta})
\ee
 Denote this
 \be\label{mapN}
 \tilde{\eta}=\mathcal{N}_F(\eta,\e)
 \ee
 
Moreover, $\mathcal{N}_F(\eta,\e)$ depends analytically on $\eta$ and $\e$ for $|\e|<\e_1$ for $\e_1$ small enough. We have $|\tilde{\eta}-1|\leq c_2|\e |$ for some $c_2>0$.
 
 Furthermore we have, in the notations (\ref{not_PHI}), (\ref{not_PSI}), (\ref{valB0}), and with
\be\label{tetp}
\theta_p=(-1)^{p+k-1}
\ee
that  for $p\geq 2$
\begin{multline}\label{serteta}
\tilde{\eta}=\mathcal{N}_F(\eta,\e)=\eta\,+\, \e^p\, \eta^{p+1}\, \frac{B_0\Phi_{p,k}(1)}{2k}\,\left(1+\theta_p\right)\\
+\, \e^{p+1}\,\eta^{p+2}\,  \frac{B_1\Phi_{p+1,k}(1)}{2k}\, \left( 1-\theta_p \right)\,+\,O(\e^{p+2})
\end{multline}
and for $p=1$ 
\begin{multline}\label{sertet1}
\tilde{\eta}=\mathcal{N}_F(\eta,\e)=\eta\,+\, \e\, \eta^{2}\, \frac{B_0\Phi_{1,k}(1)}{2k}\, (1+\theta_1)\\
+\, \e^{2}\,\eta^{3}\,  \frac{(1-\theta_1)\,B_1\,\Phi_{2,k}(1)\,+\, (1+\theta_1)/k\, B_0^2\,\Phi_{1,k}(1)^2}{2k} \,+\,O(\e^{3})
\end{multline}

\end{Lemma}

{\em{Proof of {Lemma}\,\ref{Lyp}.}}

Using (\ref{formphi}) equation (\ref{eqphi1}) is equivalent to solving the implicit equation 
\begin{multline}\label{defF}
G(\tilde{\eta},\eta,\e)=0\ \ \ {\mbox{where}}\\
G(\tilde{\eta},\eta,\e)=\tilde{\eta}^{2k}\left[ P_{J_2{F}}(1;\e\tilde{\eta})+v_{J_2{F}}(1;\e\tilde{\eta})\right]\, -\, {\eta}^{2k}\left[ P_F(1;\e\eta)+v_{F}(1;\e{\eta})\right]
\end{multline}
which is a function analytic in $(\tilde{\eta},\eta,\e)$ by Lemma\,\ref{L1}.

We have $G(\eta,\eta,0)=0$ (by (\ref{serp}), (\ref{formP0}), (\ref{servxi}))  and 
$$\frac{\partial G}{\partial \tilde{\eta}}(\eta,\eta,0)=2k\eta^{2k-1}\ne 0$$
therefore by the implicit function theorem and using the compactness of the interval $|\eta-1|\leq c_1\e_0$, equation $G(\tilde{\eta},\eta,\e)=0$ determines $\tilde{\eta}=\tilde{\eta}(\eta,\e)$ as an analytic function if  $|\e|\leq\e_1$ for $\e_1$ small enough.

Since $\tilde{\eta}(\eta,0)=\eta$ we have, for $|\e|\leq\e_1$ and $|\eta-1|\leq c_1\e_1$,
$$|\tilde{\eta}(\eta,\e)-\eta|\leq |\e|\,\sup_{|\e|\leq\e_1,|\eta-1|\leq c_1\e_1}\big|\frac{\partial\tilde{\eta}}{\partial\e}\big|=c_1'|\e| $$
therefore
$$|\tilde{\eta}(\eta,\e)-1|\leq |\tilde{\eta}(\eta,\e)-\eta|+|{\eta}-1|\leq (c_1+c_1')|\e |\equiv c_2|\e |$$
We can assume $c_2\e_1<\eta_0$ by lowering $\e_1$.

The expansion of $\tilde{\eta}(\eta,\e)$ in power series of $\e$ is found by introducing the expansions (\ref{serp}), (\ref{servxi}), (\ref{recvn}) in (\ref{defF}) and noting that $P_0(1)=1$, $P_{J_2F;2}(1)=P_2(1)$ (see (\ref{formP0}), (\ref{formPm}), (\ref{formF20})),  and that the coefficients of $B_{J_2F}$ in (\ref{valB0}) are $B_{J_2F;0}=(-1)^{p+k}B_0$ and $B_{J_2F;1}=(-1)^{p+k+1}B_1$. 

\qed

 \subsubsection{Matching at the negative $y$-axis}

Let $\eta,\tilde{\tilde{\eta}}$ satisfying (\ref{eqxyt}) and $|\eta-1|\leq c_1|\e|$ as in Lemma\,\ref{Lyp}. Let $\tilde{\eta}$ given by Lemma\,\ref{Lyp}. Consider the solution $\tilde{y_3}(x)=-\phi_{J_3F}(-x;\e,\tilde{\eta})$ as in {\bf{(iii)}}, for $x\in[-\tilde{\eta},0]$, with $\tilde{y_3}(-\tilde{\eta})=0$. Therefore  $\tilde{y_3}$ is the continuation of  $\tilde{y_2}$.

Let $\tilde{\tilde{y_4}}(x)=-\phi_{J_4F}(x;\e,\tilde{\tilde{\eta}})$ be a solution of (\ref{eqy}) as in {\bf{(iv)}}, for $x\in[0,\tilde{\tilde{\eta}}]$.

The following Lemma finds $\tilde{\tilde{\eta}}$ so that $\tilde{y_3}(0)=\tilde{\tilde{y_4}}(0)$, therefore so that $\tilde{\tilde{y_4}}$ is the continuation of $\tilde{y_3}$:

\begin{Lemma}\label{Lyp34}
Let $|\e|<\e_1$ and $|\tilde{\eta}-1|\leq c_2|\e |$ with $\e_1$ small enough so that $\tilde{\eta}$ satisfies  (\ref{eqxyt}).

There exists a unique $\tilde{\tilde{\eta}}>0$ so that 
\be\label{loweq}
-\phi_{J_4F}(0;\e,\tilde{\tilde{\eta}})=-\phi_{J_3F}(0;\e,\tilde{\eta})
\ee

Moreover, we have
\be\label{mapN3}
\tilde{\tilde{\eta}}=\mathcal{N}_{J_3F}(\tilde{\eta},\e)
\ee
where $\mathcal{N}_F$ denotes the function (\ref{mapN}) of Lemma\,\ref{Lyp}.

Therefore $\tilde{\tilde{\eta}}$ depends analytically on $\tilde{\eta}$ and $\e$ for $|\e|<\e_2$ for $\e_2$ small enough. We have $|\tilde{\tilde{\eta}}-1|\leq c_3|\e |$ for some $c_3>0$.

Furthermore, the first coefficients of the $\e$ series of $\mathcal{N}_{J_3F}(\eta,\e)$ and $\mathcal{N}_{F}(\eta,\e)$ coincide: 
\be\label{tilfirt}
\mathcal{N}_{J_3F}({\eta},\e)=\mathcal{N}_{F}({\eta},\e)+O(\e^{p+2})
\ee

\end{Lemma}

{\em{Proof.}}

Note that the equation (\ref{loweq}) for $\tilde{\tilde{\eta}}=\tilde{\tilde{\eta}}(\tilde{\eta},\e)$ is the same as the equation (\ref{defF}) for ${\tilde{\eta}}={\tilde{\eta}}({\eta},\e)$, only with $F$ replaced by $J_3F$. Noting that $P_{J_3F;2}(1)=P_2(1)$ and that $B_{J_3F;0}=\theta_pB_0$, $B_{J_3F;1}=-\theta_pB_1$  {Lemma}\,\ref{Lyp34} follows from Lemma\,\ref{Lyp}.
\qed

\subsection{The first return map}

Let $\eta$ satisfying (\ref{eta}) and $\e_2$ as in Lemma\,\ref{Lyp34}. Then $\tilde{\tilde{\eta}}$ given by Lemma\,\ref{Lyp34} is the first return to the positive $x$-axis of the solution with $x(0)=\eta,y(0)=0$ and it is analytic in $\e$ and $\eta$:
\be\label{reteta}
\tilde{\tilde{\eta}}=\mathcal{N}_{J_3F}(\mathcal{N}_{F}({\eta},\e),\e)
\ee

At this point it is enough to take $\eta=1$. Going back through the substitutions (\ref{notxy}),\,(\ref{subwu}), we obtain that the first return map to the transversal $w=z^{l+1}F(z)$ (with $z>0$ small) of the solution of (\ref{cent_foc}) with $z(0)=\e$, $w(0)=\e^{l+1}F(\e)$ is attained for 

\be\label{zret}
z\, =\, \e\,\mathcal{N}_{J_3F}(\mathcal{N}_{F}(1,\e),\e)\, \equiv\, \e\,+\,\sum_{n\geq p+1} Z_n\, \e^n
\ee
and the terms of this convergent power series in $\e$ can be calculated recursively. 

In particular, the first terms are obtained from (\ref{serteta}),\,(\ref{sertet1}),\,(\ref{tilfirt}): 
for $p\geq 2$ we have
\begin{multline}\label{zretp}
z=\e\,+\, 2\, \e^{p+1}\, \frac{B_0\Phi_{p,k}(1)}{2k}\,\left(1+\theta_p\right)\\
+\, 2\, \e^{p+2}\,  \frac{B_1\Phi_{p+1,k}(1)}{2k}\, \left( 1-\theta_p \right)\,+\,O(\e^{p+3})
\end{multline}
and for $p=1$ 
\begin{multline}\label{zret1}
z=\e\,+\, 2\, \e^{2}\, \frac{B_0\Phi_{1,k}(1)}{2k}\, (1+\theta_1)\\
+\, 2\,\,\e^{3}\,  \left[\, \frac{(1-\theta_1)\,B_1\,\Phi_{2,k}(1)\,+\, (1+\theta_1)/k\, B_0^2\,\Phi_{1,k}(1)^2}{2k}\, \right.\\
+\left.\, \left(\frac{B_0\Phi_{1,k}(1)}{2k}\, (1+\theta_1)\ \right)^2\right] \,+\,O(\e^{3})
\end{multline}

\subsection{Closed trajectories}\label{center}

Relation (\ref{zret}) together with the constructive method presented allows to calculate, recursively, the coefficients $Z_n$, and therefore to decide whether the origin is a center or a focus: the origin is a center if and only if all $Z_n$ vanish.

As a practical matter, for each value of $p=l-k+1$ all the power series should first be properly ordered (the order in which the terms appear does depend on the value of $p$) and then the series can be calculated (at least in principle) to any order $n_0$. The conditions $Z_n=0$ for $n\leq n_0$ can be written in terms of the function $F$, and they are necessary conditions for the fixed point to be a center.

The first such conditions follow from the first two nontrivial terms calculated here: using (\ref{zretp}),\,(\ref{zret1}) for the origin to be a center we must have $Z_{p+1}=Z_{p+2}=0$ which implies $B_0\left(1+\theta_p\right)=0$ and $B_1\left( 1-\theta_p \right)=0$, which in turn means that if $p+k-1$ is odd then we must have $F'(0)=0$, while if $p+k-1$ is even, then we must have $F(0)=0$.

\ 

\ 

{\bf{Acknowledgements.}} The author is grateful to Chris Miller for suggesting the problem, to Barbara Keyfitz for very interesting discussions, and to Christiane Rousseau for illuminating e-mail correspondence.

\end{document}